# BOUNDED RANK SUBGROUPS OF COXETER GROUPS, ARTIN GROUPS AND ONE-RELATOR GROUPS WITH TORSION

ILYA KAPOVICH AND PAUL SCHUPP

ABSTRACT. We obtain a number of results regarding freeness, quasiconvexity and separability for subgroups of Coxeter groups, Artin groups and one-relator groups with torsion.

## 1. INTRODUCTION

Showing that various subgroups of given groups are free is a central theme of combinatorial and geometric group theory. Two examples which immediately come to mind are the fact that the Kurosh Subgroup Theorem shows that a subgroup of a free product which has trivial intersection with all conjugates of the factors is free, and the celebrated Freiheitssatz of Magnus showing that if $H$ is a subgroup of a one-relator group which is generated by a subset of the generators omitting a generator occurring in the cyclically reduced relator then $H$ is free. Indeed, Arzhantseva and Olshanskii [6] showed that, generically, for most groups given by $m$ generators and $n$ defining relators, any subgroup with fewer than $m$ generators is free. Their proof introduces a minimization argument on subgroup graphs. In subsequent work Arzhantseva [1, 2, 3] applied this approach to prove a number of other results about "generic" properties of finitely presented groups (see also the related work of Bumagina [9]).

We believe that the Arzhantseva-Olshanskii method deserves to be much more widely known. In this paper we show that this approach can be combined with other techniques to yield precise (rather than probabilistic) results about subgroups of bounded rank in some well-known classes of groups, namely, Coxeter groups, Artin groups and one-relator groups with torsion.

A subgroup $H$ of a word-hyperbolic group $G$ is *quasiconvex* for some (and hence for any) finite generating set $A$ of $G$ and the corresponding word-metric there is a fixed bound $B$ such that for any geodesic representative of any element of $H$, the group element represented by any initial segment is within distance $B$ of an element of $H$. A word-hyperbolic group $G$ is *locally quasiconvex* if all its finitely generated subgroups are quasiconvex. Quasiconvexity in the context of Kleinian groups essentially corresponds to geometric finiteness. More precisely, if $G$ is a geometrically finite group of isometries of $\mathbb{H}^n$ without parabolic elements then a subgroup $H \leq G$ is quasiconvex if and only if $H$ is geometrically finite [57]. Free group of finite rank, hyperbolic surface groups, many 3-manifold groups as well as many small cancellation groups are known to be locally quasiconvex. Locally quasiconvex groups enjoy a number of particularly good algebraic, geometric and algorithmic properties and their study plays an important role in the theory of word-hyperbolic groups [44, 34, 35, 22, 23, 24, 58].

By a $k$-generated group we mean a group which can be generated by $k$ or fewer elements.

**Theorem A.** *Let $k$ be a positive integer and let*

$$(\dagger) \qquad G = \langle a_1, \ldots, a_n \mid a_i^2 = 1, i = 1, \ldots, n; (a_i a_j)^{m_{ij}} = 1, 1 \leq i < j \leq n \rangle$$

*be a Coxeter group such that for all $i < j$ we have $m_{ij} \geq 3k+1$. Then the following holds:*

1. *Every torsion-free $k$-generated subgroup of $G$ is free and quasiconvex in $G$.*
2. *Every $k$-generated subgroup of $G$, which does not contain any conjugates of the generators $a_i$, is quasiconvex in $G$.*

3. *If all finite $m_{ij}$ are even then every subgroup of $G$ generated by at most $(k+1)/2$ elements, is quasiconvex in $G$.*

Note that in the presentation above we allow some $m_{ij} = \infty$, which means that there is no relation between $a_i$ and $a_j$. All Coxeter groups and Artin groups which we consider will be of *extra-large type*, that is, all $m_{ij}$ are at least 4.

Recall that a subgroup $H \leq G$ is said to be *separable* if for any element $g \in G - H$ there exists a homomorphism $\psi : G \to K$ onto a finite group $K$ such that $\psi(g) \notin \psi(H)$, that is, if $H$ equals the intersection of all finite index subgroups of $G$ which contain $H$. A group is said to be *subgroup separable* if all its finitely generated subgroups are separable. (Some authors use the term LERF for this concept.) We obtain quasiconvexity in the theorem above and separability in the theorem below because the Arzhantseva-Olshanskii minimization technique produces a subgroup graph which, as explained later, has a finite 2-completion by the results of Schupp [52]. For separability we also need the following

**Definition 1.1.** [The Separability Condition] $G$ is of extra-large type, all finite $m_{ij}$ in the presentation (†) are even and, whenever three distinct generators $a_i, a_j, a_k$ are all pair-wise related in (†) (that is $m_{ij} < \infty$, $m_{ik} < \infty$ and $m_{jk} < \infty$), then $m_{ij}$ is also divisible by 3 (and thus divisible by 6).

**Theorem B.** *Let $k \geq 2$ and $G$ be as in (†) such that for each $i < j$ we have $m_{ij} \geq 3k+7$ and $G$ satisfies the Separability Condition.*

*Then*

1. *Every $k$-generated subgroup of $G$ which does not contain any conjugates of the generators $a_i$ is separable in $G$.*
2. *Every subgroup of $G$ generated by at most $(k+1)/2$ elements is separable in $G$.*

Since a group $G$ in both Theorem A and Theorem B is word-hyperbolic, all possible definitions of quasiconvexity for subgroups of $G$ coincide. In particular a subgroup $H \leq G$ is quasiconvex if and only if $H$ is finitely generated and quasi-isometrically embedded in $G$ and if and only if $H$ is rational with respect to an automatic language for $G$ (see [4, 28, 33, 43]). Thus the membership problem for a quasiconvex subgroup of $G$ is solvable in linear time. Quasiconvex subgroups of $G$ are finitely presentable and themselves word-hyperbolic. Also, the intersection of finitely many quasiconvex subgroups of $G$ is again quasiconvex.

Separability with respect to quasiconvex subgroups is related to separability of geometrically finite subgroups for Kleinian groups. The latter property plays a prominent role in 3-dimensional topology and has been the subject of extensive research [5, 23, 37, 38, 39, 46, 53, 60]. Long and Reid [39] recently obtained some results regarding separability with respect to geometrically finite subgroups of Coxeter groups generated by reflections in the sides of a non-compact hyperbolic simplex of finite volume. Some results related to subgroup separability in right-angled Coxeter groups are also obtained by Gitik in [25].

Schupp [52] used the perimeter reduction idea of McCammond and Wise [44] to define an extensive class of Coxeter groups of extra-large type for which all finitely generated subgroups have subgroup graphs admitting a finite 2-completion. Therefore these groups are locally quasiconvex and those which also satisfy the Separability Condition are subgroup separable.

The Arzhantseva-Olshanskii minimization technique is complementary to perimeter reduction and the hypotheses required are different. Perimeter reduction requires a *distributive* small cancellation condition which involves how each generator is distributed among *all* the defining relators. Thus how large an $m_{i,j}$ has to be depends on how many relators involve the generators $a_i$ and $a_j$. If perimeter reduction works it gives results on *all* finitely generated subgroups. The minimization technique requires the exponents $m_{i,j}$ to increase with the rank $k$ of the subgroup but is completely independent of how many relators involve given generators.

Appel and Schupp [7] introduced the use of small cancellation methods into the study of Artin groups. Among their results they proved that Artin groups of extra-large type are torsion-free and that the squares of the generators freely generate a free subgroup. This is the source of Tits' Conjecture, recently proven by Crisp and Paris [17], that in any Artin group a relation between the squares of the generators is a consequence solely of the relators which are commutators. Perimeter reduction does not seem applicable to Artin groups, but one *can* use the Arzhantseva-Olshanskii technique together with the appropriate small cancellation theory.

**Theorem C.** *Let $k \geq 2$ be an integer. Let*

$$(1) \qquad G = \langle a_1, \ldots, a_n \,|\, u_{ij} = u_{ji}, \text{ where } 1 \leq i < j \leq n, \rangle$$

*be an Artin group where for $i \neq j$*

$$u_{ij} := \underbrace{a_i a_j a_i \ldots}_{m_{ij} \text{ terms}}$$

*and where $m_{ij} = m_{ji} \geq 5k$ for each $i < j$.*

*Suppose $H \leq G$ is a $k$-generated subgroup such that $H$ has trivial intersection with all conjugates of every two-generator subgroup $G_{ij} = \langle a_i, a_j \rangle$ for which $m_{ij} < \infty$. Then $H$ is free.*

Few general results seem to be known about the subgroup structure of Artin groups. Most subgroup theorems for Artin groups deal with subgroups generated by subsets of the generators or with subgroups of finite index or with kernels of homomorphisms to abelian groups (see for example [8, 10, 11, 12, 13, 16, 17, 18, 20, 42, 47]). Hsu and Wise [31] recently obtained results about the separability of quasiconvex subgroups for a restricted class of right-angled Artin groups.

The consideration of free subgroups is, of course, completely crucial to the study of one-relator groups. The class of one-relator groups with torsion is often the most tractable class of one-relator groups because of the power of Newman's Spelling Theorem [45]. Steve Pride [48] proved that every torsion-free two-generator subgroup of a one-relator group with torsion is free.

**Definition 1.2.** If $r$ is a cyclically reduced word of a free group $F(a_1, \ldots, a_n)$, the *letter bound* $b(r)$ is the maximum number of times any letter $a \in \{a_1, \ldots, a_n, a_1^{-1}, \ldots, a_n^{-1}\}$ occurs in $r$. (Given a letter, we do not count the occurrences of its inverse.).

For example, $b(a_1 a_2 a_1^{-1} a_2^{-1}) = 1$ and $b(a_1^2 a_2^d) = |d|$ provided $|d| \geq 2$.

**Theorem D.** *Let $k \geq 2$ be an integer and let*

$$(2) \qquad G = \langle a_1, \ldots, a_n \,|\, r^m = 1 \rangle$$

*be a one-relator group where, as usual, $r$ is a nontrivial cyclically reduced word which is not a proper power. Let $b = b(r)$. If $m \geq b(6k - 2) + 2$ then every $k$-generated torsion-free subgroup of $G$ is free.*

The results of McCammond and Wise [44] and Hruska and Wise [30] using perimeter reduction imply that many of the one-relator groups with torsion satisfying conditions of Theorem D are locally quasi-convex. But again the hypotheses are not the same. Since there are arbitrarily long cyclically reduced words $r$ having $b(r)$ equal to any given positive value, there are many groups covered by Theorem D to which perimeter reduction does not apply.

The main idea of the minimization technique can be loosely described as follows. For a group $G$ with a fixed finite generating set $A$ subgroups of $G$ can be represented by graphs whose edges are labeled by the elements of $A$. A $k$-generated subgroup of $G$ can always be represented by a graph whose Euler characteristics is $\geq 1 - k$. Standard Stallings folds as well as additional Arzhantseva-Olshanskii moves preserve the subgroup represented by a graph. Given a $k$-generated subgroup $H \leq G$ we choose among the labeled graphs of Euler characteristic $\geq k$ representing $H$ a graph with minimal complexity (which needs to be appropriately defined, depending on the group $G$). We then study this minimal graph $\Gamma$ to obtain the desired results. In many instances it is possible to use minimality together with the appropriate small cancellation theory to show that the canonical epimorphism from $\pi_1(\Gamma)$ to $H$ is a monomorphism and hence $H$ is free.

## 2. Representing subgroups by labeled graphs

**Convention 2.1.** Let $A = \{a_1, \ldots, a_n\}$ be a finite alphabet, which will be the set of generators of the group $G$ under consideration. In discussing Artin groups and one-relator groups with torsion, $F$ denotes the free group $F(A) = F(a_1, \ldots, a_n)$ on $A$. In our discussion of Coxeter groups, $F$ denotes the free product of cyclic groups of order two:

$$\bar{F}(a_1, \ldots, a_n) = \bar{F}(A) := \bigstar_{i=1}^{n} \langle a_i \mid a_i^2 = 1 \rangle.$$

In the context of Coxeter groups a word $w$ will be called *reduced* if it does not contain subwords of the form $a_i a_i$. In the context of quotients of $F(A)$ a word $w$ will be called *reduced* if it freely reduced, that is if $w$ does not contain subwords of the form $a_i a_i^{-1}$ or $a_i^{-1} a_i$. A reduced word $w$ is said to be *cyclically reduced* if all cyclic permutations of $w$ are reduced. The number of letters in a word $w$ will be called the *length* of $w$ and denoted $|w|$. The element of $G$ represented by a word $w$ will be denoted $\overline{w}$ (where $G$ is either a Coxeter group or an Artin group or a one-relator group with torsion, depending on the context). We will usually use $w_*$ and $w^*$ to denote a cyclic permutation of the word $w$. If $F = F(A)$ and $w$ is a word, the $w^{-1}$ will denote the standard formal inverse of $w$. If $w = y_1 \ldots y_l$ is a word in $\bar{F}(A)$, where each $y_h \in A$, we set its formal inverse to be $w^{-1} := y_l \ldots y_1$.

Following the approach of Stallings [55], we use labeled graphs to study finitely generated subgroups of quotients of either a free group or a free product of cyclic groups of order two [55, 52, 6, 32].

**Definition 2.2.** An *$F$-graph* $\Gamma$ consists of an underlying oriented graph where every edge $e$ is labeled by a nontrivial reduced word $\mu(e)$ in such a way that $\mu(e^{-1}) = \mu(e)^{-1}$ for every edge $e$ of $\Gamma$. We allow multiple edges between vertices as well as edges which are loops.

An *$F$-graph* $\Gamma$ is said to be *non-folded* if there exists a vertex $v$ and two distinct edges $e_1, e_2$ with origin $v$ such that the words $\mu(e_1)$ and $\mu(e_2)$ have a nontrivial common initial segment (in the Coxeter case we also require that if $|\mu(e_1)| = |\mu(e_2)| = 1$ then $e_1 \neq e_2^{-1}$.) Otherwise $\Gamma$ is said to be *folded*.

Every edge-path $p$ in $\Gamma$ has a label which is a word in $A$. We shall denote this label by $\mu(p)$. The number of edges in $p$ will be called the *length* of $p$ and denoted $|p|$. For $F = F(A)$ a path $p$ in an $F$-graph $\Gamma$ is said to be *reduced* if it does not contains subpaths of the form $e, e^{-1}$ where $e$ is an edge of $\Gamma$. For $F = \bar{F}(A)$ a path $p$ in an $F$-graph $\Gamma$ is said to be *reduced* if it does not contains subpaths of the form $e, e^{-1}$ and if $p$ does not contain subpaths of the form $e, e$ where $e$ is an edge of $\Gamma$ with $|\mu(e)| = 1$.

The following statement is obvious:

**Lemma 2.3.** *Let $\Gamma$ be an $F$-graph. Then $\Gamma$ is folded if and only if the label of any reduced path in $\Gamma$ is a reduced word.*

We need the following straightforward modification of Stallings' folding moves:

**Definition 2.4** (Fold)**.** Let $\Gamma$ be an $F$-graph. Suppose $e_1 \neq e_2$ are distinct edges of $\Gamma$ with labels $w_1$ and $w_2$ accordingly and with a nontrivial common initial vertex $v$ (If $F = \bar{F}(A)$ and $|\mu(e_1)| = |\mu(e_2)| = 1$ we also suppose $e_1 \neq e_2^{-1}$). Suppose $x$ is the maximal common initial segment of $w_1$ and $w_2$ and assume that $|x| > 0$. Thus we have graphic equalities $w_1 = x u_1$ and $w_2 = x u_2$ where the product $u_2^{-1} u_1$ is reduced. For $i = 1, 2$ if $|u_i| > 0$ we subdivide $e_i$ into two consecutive edges labeled $x$ and $u_i$ accordingly. After that we fold the two edges labeled $x$ originating at $v$ into a single edge $e$ labeled $x$.

The resulting $F$-graph $\Gamma'$ is said to be obtained from $\Gamma$ by a *fold*.

The following statement immediately follows from the definitions, exactly as in [55]:

**Proposition 2.5.** *Let $\Gamma$ be a connected graph and suppose that $\Gamma'$ is obtained from $\Gamma$ by a fold. Then the Euler characteristic of $\Gamma'$ is no less than that of $\Gamma$.*

In addition to the Stallings folding moves, we need the following two transformations of labeled graphs introduced by Olshanskii and Arzhantseva [6]. Recall that we are working with a fixed presentation of a quotient $G$ of $F$ so "relator" is defined in terms of this given presentation.

**Definition 2.6** (Completing a relator cycle: move $A0$). Let $p$ be path in an $F$-graph $\Gamma$ with an initial vertex $x$, a terminal vertex $y$ and label $\mu(p) = v$. Suppose $v'$ is a reduced word such that $\overline{v} = \overline{v'} \in G$. We modify $\Gamma$ by attaching a new edge going from $x$ to $y$ labeled by the word $v'$.

**Definition 2.7** (Removing a simple edge from a relator cycle: move $A1$). Let $e$ be an edge in a labeled $F$-graph $\Gamma$ with an initial vertex $x$, a terminal vertex $y$ and label $\mu(e) = v$. Suppose there exists a path $p'$ in $\Gamma$ from $x$ to $y$ with label $\mu(p') = v'$ such that $p'$ does not contain $e$ or $e^{-1}$ and such that $\overline{v} = \overline{v'} \in G$.

We modify $\Gamma$ by removing the edge $e$ while keeping the vertices $x$ and $y$.

Note that $A0$ decreases the Euler characteristic by one and $A1$ increases the Euler characteristic by one.

We also need the following two moves:

**Definition 2.8** (Subdividing an edge: move $A2$). Let $e$ be an edge in a labeled $F$-graph $\Gamma$ with an initial vertex $x$, a terminal vertex $y$ and label $\mu(e) = u$.

Suppose $u = u_1 u_2$ is a graphical equality where both $u_1, u_2$ are nontrivial reduced words (thus the product $u_1 u_2$ is reduced).

We modify $\Gamma$ as follows: Add a new vertex $z$ subdividing the edge $e$ into two new edges $e_1$ and $e_2$ labeled $u_1$ and $u_2$ so that $e_1$ goes from $x$ to $z$ and $e_2$ goes from $z$ to $y$.

**Definition 2.9** (Removing a vertex of degree two: move $A3$). Suppose $\Gamma$ is an $F$-graph and $z$ is a vertex of degree two in $\Gamma$ such that no loop-edge is based at $z$. Let $e_1$ and $e_2 \neq e_1^{-1}$ be the two distinct edges with origin $z$ and labels $u_1^{-1}$ and $u_2$ accordingly. (Note that $e_i$ is not a loop.) Suppose that $u_1 u_2 \neq 1$ in $F$. Let $u$ be the reduced form of $u_1 u_2$, so that $|u| > 0$. Denote the terminal vertex of $e_1$ by $x$ and the terminal vertex of $e_2$ by $y$.

We modify $\Gamma$ as follows: remove the vertex $z$, combining the edges $e_1, e_2$ into a single new edge $e$ with label $u$ going from $x$ to $y$.

Exactly as in the case of the standard Stallings graph, where edges are labeled by individual letters from $A$ rather than by words, we obtain the following:

**Proposition-Definition 2.10.** Let $\Gamma$ be a connected $F$-graph with a base-vertex $x_0$. Then the labeling of paths gives rise to a homomorphism

$$\phi : \pi_1(\Gamma, x_0) \to G$$

such that for every path $p$ from $x_0$ to $x_0$ we have $\phi([p]) = \overline{\mu(p)}$ (where $[p]$ stands for the equivalence class of $p$ in $\pi_1(\Gamma, x_0)$). In this case we will say that $H = \phi(\pi_1(\Gamma, x_0)) \leq G$ is the *subgroup represented by* $(\Gamma, x_0)$.

Moreover, the following holds:

1. If $\Gamma$ is finite then $image(\phi)$ is finitely generated. Moreover if $\Gamma$ has Euler characteristic $1 - k$, then $\pi_1(\Gamma, x_0)$ and hence $\phi(\pi_1(\Gamma, x_0))$ can be generated by $k$ elements.
2. Every finitely generated subgroup of $G$ can be represented in this fashion for some finite connected $\Gamma$. Moreover, if $H \leq G$ is $k$-generated, then $H$ can be represented by a connected $F$-graph of Euler characteristic $\geq 1 - k$.
3. If $x_1$ is another vertex of $\Gamma$ then the pairs $(\Gamma, x_0)$ and $(\Gamma, x_1)$ define conjugate subgroups of $G$.

The following simple fact plays an important role in our approach.

**Proposition 2.11.** *Let $\Gamma$ be an $F$-graph with a base-vertex $x_0$. Suppose $\Gamma'$ is obtained from $\Gamma$ by a finite sequence of folds and moves $A0, A1, A2, A3$, so that $x_0'$ is the image of $x_0$. Then the pairs $(\Gamma, x_0)$ and $(\Gamma', x_0')$ define the same subgroup of $G$.*

*Proof.* It suffices to prove the result for a single fold or move. For a fold the statement follows exactly as in the case of the standard Stallings folds [55, 32]. For the moves $A0, A1$ the proof is exactly the same as for Lemma 1 in [6]. The statement regarding the moves $A2, A3$ is an easy corollary of the definitions. $\square$

We need a measure of complexity for $F$-graphs. When considering Coxeter groups and one-relator groups with torsion, we use the following measure. (We shall use a slightly more complicated measure when studying Artin groups.)

**Definition 2.12** (Complexity). Let $\Gamma$ be an $F$-graph. We partition the set $E(\Gamma)$ of *oriented* edges of $\Gamma$ as $E(\Gamma) = E^+ \cup E^-$ where for any $e \in E(\Gamma)$ we have $e \in E^+ \iff e^{-1} \in E^-$ and where $E^+ \cap E^- = \emptyset$.

Let $l = \sum_{e \in E^+} |\mu(e)|$ and let $q$ be the number of vertices in $\Gamma$. We define the *complexity* $\sigma(\Gamma)$ to be $\sigma(\Gamma) = (l, q)$.

We order complexities lexicographically. That is $(l, q) < (l', q')$ if one of the following occurs:

1. $l < l'$; or
2. $l = l'$ and $q < q'$.

We need the following simple graph-theoretic lemma.

**Lemma 2.13.** *Let $Y$ be a finite connected non-oriented graph with Euler characteristic $1 - k$, where $k \geq 2$ such that each vertex of $Y$ has degree at least three. Then*

1. *A maximal subtree of $Y$ contains at most $2k - 3$ edges.*
2. *If $p$ is a simple path in $Y$ then $p$ contains at most $2k - 3$ edges.*

*Proof.* Since any simple path can be extended to a maximal subtree $T$, the second statement immediately from the first. To show that a maximal subtree has at most $2k - 3$ edges, let $T$ be a maximal subtree in $Y$, say with $x$ edges. By assumption $Y - T$ has $k$ edges, so that $Y$ has $E := k + x$ edges. Note that $T$ has $x + 1$ vertices. Since $T$ is a maximal subtree of $Y$, it follows that every vertex of $Y$ is present in $T$ and so $Y$ has $V := x + 1$ vertices. Since the degree of every vertex is at least three, we have $2E \geq 3V$, that is $2(k + x) \geq 3(x + 1)$. This implies $x \leq 2k - 3$, as required. $\qquad\square$

## 3. SMALL CANCELLATION THEORY AND COXETER GROUPS

We need only a few basic facts from small cancellation theory and how they apply to Coxeter groups as well as Artin groups. Small cancellation theory over a free product of cyclic groups of order two is essentially the same as the usual theory over a free group. We follow Lyndon-Schupp [40] and Appel-Schupp [7] and the reader is referred to these sources for the background information. A set $R$ of reduced words in $F$ (where $F$ is either $\bar{F}(A)$ or $F(A)$) is *symmetrized* if $R$ is closed under taking inverses and cyclic permutations. A *piece* (relative to the set $R$) is a subword $u$ such that there are *distinct* elements

$$r = ur_1' \text{ and } r_2 = ur_2'$$

in $R$. The set $R$ satisfies the *condition* $C(p)$ if no element of $R$ is a product of fewer than $p$ pieces. Also, $R$ ratifies the *condition* $C'(\lambda)$ if the length of any piece $u$ is less than $\lambda$ times the length of any relation in which $u$ occurs as a subword.

As stated in the introduction, we view the Coxeter presentation (†) as a presentation over the free product of cyclic groups of order two $F = \bar{F}(A) = *_{i=1}^n G_i$ where $G_i = \langle a_i \mid a_i^2 = 1 \rangle$, with defining relations $r_{ij} := (a_i a_j)^{m_{ij}}$:

$$(\ddagger) \qquad\qquad G = F/\langle r_{ij} = 1, i < j \rangle.$$

The symmetrized set for a Coxeter presentation just consists of the words $(a_i a_j)^{m_{ij}}$ and $(a_j a_i)^{m_{ij}}$, where $m_{ij} < \infty$ and $i < j$. Thus all pieces for presentation (‡) have length one and are individual letters of $A$. All the Coxeter groups which we consider are of extra-large type, that is, all $m_{ij} \geq 4$. Such a presentation therefore satisfies $C(8)$ and $C'(1/7)$. The "basic fact" of small cancellation theory is the following [40]:

**Proposition 3.1.** *Let $G = <A; R>$ where $R$ is symmetrized and satisfies $C(6)$. Then:*

1. *If $w$ is a nontrivial word which is equal to the identity in $G$ then $w$ contains a subword of an element of $R$ which is a relator with at most three pieces missing.*
2. *If $R$ satisfies $C'(1/6)$ then Dehn's Algorithm solves the word problem for $G$. Thus if $R$ is finite then $G$ is word-hyperbolic.*

**Definition 3.2.** Let $w$ be a reduced word. We say that $w$ is a *weakly Dehn-reduced* word if $w$ does not contain a subword $v$ such that $v$ is also a subword of some $r_{ij}^*$ with $|v| \geq |r_{ij}| - 3$.

In a Coxeter group of extra-large type, a nontrivial word $w$ which is weakly Dehn-reduced is not equal to the identity. Also, the Torsion Theorem for $C'(1/6)$ groups gives a simple explicit description of torsion elements [7]:

**Proposition 3.3.** *Let $G$ be given by presentation* (†) *and suppose that all $m_{ij} \geq 4$. Then each $a_i$ has order two in $G$ and each $a_i a_j$ has order $m_{ij}$ in $G$. Moreover, every nontrivial element of finite order in $G$ is conjugate to either some $a_i$ or to some $(a_i a_j)^d$, where $0 < |d| < m_{ij} < \infty$.*

Note that this characterization shows that a subgroup $H$ of $G$ is torsion-free if and only if it has trivial intersection with all conjugates of the cyclic subgroups on the generators and with all conjugates of the subgroups $< a_i, a_j >$ for which $m_{ij} < \infty$.

Until the end of Section 4, unless specified otherwise, we assume that $G$ is given by presentation (‡) where for all $i < j$ we have $m_{ij} \geq 4$. Thus Proposition 3.1 applies and $G$ is word-hyperbolic. Denote by $X = X(G, A)$ the Cayley graph of $G$ with respect to $A$. That is to say, the vertex set of $X$ is $G$ and whenever $g_1, g_2 \in G$ are such that $g_1 =_G g_2 a_i$ (and hence $g_2 =_G g_1 a_i$), the vertices $g_1$ and $g_2$ are connected by a *single* edge labeled $a_i$ in $X$.

We also denote by $d_A$ the word-metric on $X$ corresponding to $A$. For an element $g \in G$ denote $|g|_A := d_A(1, g)$. A path $p$ in $X(G, A)$ is said to be *geodesic* if the length of $p$ is equal to the distance between its endpoints. A word $w$ is said to be *geodesic* if $|w| = |\overline{w}|_A$. Thus the label of an edge-path is a geodesic word if and only if the path is geodesic. These notations and conventions will also be fixed till the end of Section 4.

## 4. The main results about Coxeter groups

We need the following definition from Schupp [52].

**Definition 4.1** (The Relator Path Property). Suppose $\Gamma$ is an $F$-graph for a Coxeter group $G$ of extra-large type. We say that a closed edge-path $p$ with label $r_{ij}^*$ is a *relator cycle* in $\Gamma$. The graph $\Gamma$ is said to have the *Relator Path Property*, if $\Gamma$ is folded and the following holds.

Let $\Gamma'$ be the graph obtained from $\Gamma$ by performing edge-subdivisions so that the label of every edge in $\Gamma'$ has length one. Then whenever $p$ is a path in $\Gamma'$ whose label $\mu(p)$ is a subword of some $r_{ij}^*$ missing at most three letters, then $p$ is a subpath of a relator-cycle in $\Gamma'$.

Note that if $\Gamma$ is folded and the label of every reduced edge-path in $\Gamma$ is weakly Dehn-reduced, then $\Gamma$ clearly has the Relator Path Property. We refer to [52] for a precise definition of a graph being *2-complete* for presentation (‡) but this means that any two successive edges labeled respectively by $a_i$ and $a_j$ where $m_{ij} \leq \infty$ are part of a relator cycle except possibly in one very special case. We recall the following result of [52].

**Proposition 4.2.** *Let $G$ be given by presentation* (†) *where all $m_{ij} \geq 4$ (and hence $G$ is word-hyperbolic). Suppose $\Gamma$ is a connected folded finite $F$-graph with the Relator Path Property. Let $x_0 \in \Gamma$ be a base-vertex and let $H \leq G$ be the subgroup defined by the pair $(\Gamma, x_0)$. Then:*

1. *The graph $\Gamma$ admits a finite $2$-completion, that is, a finite folded connected $2$-complete graph $\Delta_2(\Gamma)$ which contains as a subgraph the graph obtained from $\Gamma$ by subdividing every edge of $\Gamma$ into edges with labels of length one.*
2. *The subgroup $H$ is quasiconvex in $G$.*
3. *If in addition* (†) *satisfies the Separability Condition, then there is a homomorphism $\psi$ from $G$ into the symmetric group on the vertices of $\Delta_2$ with the following property. If $x$ is a vertex of $\Delta_2$ and $g \in G$ is an element represented by the label of some path from $x_0$ to $x$ in $\Delta_2$ then the permutation $\psi(g)$ takes $x_0$ to $x$. In particular, $\psi(H)$ fixes $x_0$ and if $x \neq x_0$ then $\psi(g) \notin \psi(H)$.*

Our main technical tool is the following:

**Proposition 4.3.** *Let $k \geq 2$ and $G$ be as in (†) such that for each $i < j$ we have $m_{ij} \geq 3k+1$. Suppose $H \leq G$ is a $k$-generated nontrivial subgroup (and hence $H$ can be represented by a graph with Euler characteristic $\geq 1 - k$).*

*Among all $F$-graphs with Euler characteristic $\geq 1 - k$ representing subgroups conjugate to $H$ in $G$ choose a graph $\Gamma$ with the smallest possible complexity.*

*Then $\Gamma$ is folded, connected, has no degree-one vertices. Also, either $\Gamma$ consists of a single loop or $\Gamma$ has no vertices of degree two.*

*Moreover:*

*(a) If $H$ is torsion-free then the label of every reduced edge-path in $\Gamma$ is weakly Dehn-reduced and hence $\Gamma$ has the Relator Path Property.*

*(b) If $H$ contains no conjugates of the generators $a_i$ then $\Gamma$ has the Relator Path Property.*

*Proof.* By the minimal choice of $\Gamma$ it is obvious that $\Gamma$ is folded, connected and has no degree-one vertices. In particular, this means that the label of any reduced path in $\Gamma$ is a reduced word. Similarly, the minimality assumption implies that $\Gamma$ does not have any vertices of degree two unless $\Gamma$ is a single edge which is a loop.

We will first establish part (a). Suppose that $H$ is torsion-free but there exists a nontrivial reduced edge-path in $\Gamma$ whose label is not weakly Dehn-reduced.

Then there exists a reduced path $p$ in $\Gamma$ with $\mu(p) = uvu'$, where $v$ is also a subword of some $r_{ij}^*$ with $|r_{ij}| - 3 \leq |v| \leq |r_{ij}| - 1$. By performing subdivision moves $A2$ on the edges $e, e'$ if necessary, we obtain a graph $\Gamma'$ such that there is a reduced edge-path $p'$ with $\mu(p') = v$. Note that if $\sigma(\Gamma) = (l, q)$ and $\sigma(\Gamma') = (l', q')$ then $l = l'$ and $q' \leq q + 2$. Also, $\Gamma'$ is folded since $\Gamma$ was folded and the Euler characteristics of $\Gamma$ and $\Gamma'$ are equal.

**Claim 1.** The path $p'$ is a simple path in $\Gamma'$.

Indeed, suppose this is not the case. Recall that since $p'$ is a reduced path we can choose a nonempty subpath $p_1$ of $p'$ such that $p_1$ is a simple closed path in $\Gamma'$.

Suppose first that the length of $|\mu(p_1)| = 1$. Then $\Gamma'$ possesses an edge-loop with label $a_i$ or $a_j$. This implies that $H$ contains an element conjugate to either $a_i$ or $a_j$ contradicting our assumption that $H$ is torsion-free. Suppose now that $|\mu(p_1)| > 1$. If the first and the last letters of $\mu(p_1)$ are the same, then $\Gamma'$ is not folded, contrary to our assumptions. Thus the first and the last letters of $\mu(p_1)$ are different. Without loss of generality we may assume $\mu(p_1)$ starts with $a_i$ and ends with $a_j$. Since $v$ is a subword of $r_{ij}^*$ with $|v| \leq |r_{ij}| - 1$, we conclude that $\mu(p_1) = (a_i a_j)^d$ where $1 \leq d < m_{ij}$. However $p_1$ is a closed path in $\Gamma'$, and therefore $H$ contains an element conjugate to $(a_i a_j)^d$. Again, this contradicts our assumption that $H$ is torsion-free. Thus Claim 1 is verified.

We now know that $p'$ is a simple path in $\Gamma'$. Suppose first that $\Gamma'$ has negative Euler characteristic $1 - k_0 < 0$, so that $k_0 \geq 2$ and the fundamental group of $\Gamma'$ is free of rank $k_0$. Recall that in $\Gamma$ every vertex has degree at least three and that $\Gamma'$ was obtained from $\Gamma$ by subdividing at most two edges. Hence Lemma 2.13 implies that $|p'| \leq 2k - 1$, that is $p'$ consists of at most $(2k - 3) + 2 = 2k - 1$ edges.

Therefore there is an edge $f$ of $p'$ such that $|\mu(f)| \geq |\mu(p')|/(2k-1)$. Recall that $|p'| = |v| \geq 2m_{ij} - 3$. Hence

$$|\mu(f)| \geq \frac{|\mu(p')|}{2k-1} \geq \frac{2m_{ij} - 3}{2k-1} > 3,$$

where the last inequality holds since by assumption $m_{ij} \geq 3k+1$. Denote the first vertex of $p'$ by $x$ and the last vertex of $p'$ by $y$.

We now perform a move of type $A0$ by adding an edge from $x$ to $y$ with label of length at most three corresponding to the part of $r_{ij}^*$ missing in $v$. Next we perform a type-$A1$ move and remove the edge $f$ with label of length at least four. By Proposition 2.10 the resulting graph $\Gamma_1$ represents a subgroup conjugate to $H$ in $G$. Also $\Gamma, \Gamma'$ and $\Gamma_1$ have equal Euler characteristics. However, $\Gamma_1$ has smaller complexity than does $\Gamma$, which contradicts the choice of $\Gamma$. Thus the statement of the proposition has been verified in this case.

Suppose now that $\Gamma$ and $\Gamma'$ have zero Euler characteristic and infinite cyclic fundamental group. Since $\Gamma$ has no degree-one vertices, this implies that $\Gamma$ is a single loop-edge. Thus $\Gamma'$ is a topological circle subdivided in at most three edges. Therefore $p'$ is a non-loop edge in $\Gamma'$. Again, let $x$ and $y$ be the first and the last vertices of $p'$ accordingly. Recall that $p'$ is labeled by a subword $v$ of $r_{ij}^*$ missing at most three letters and with

$$|v| \geq 2m_{ij} - 3 \geq 2(3k+1) - 3 = 6k - 1 > 3.$$

As in the previous case, we perform move $A0$ and attach an edge from $x$ to $y$ with label of length at most three corresponding to the missing in $v$ portion of $r_{ij}^*$. Then we perform move $A1$ and remove the edge $p'$. The resulting graph represents a subgroup conjugate to $H$ but has smaller complexity than $\Gamma$, yielding a contradiction. Thus part (a) is proved.

The prove of part (b) is very similar and we will briefly indicate where the changes in the argument need to be made. Suppose $H$ does not contain any conjugates of the generators $a_i$ but that $\Gamma$ does not have the Relator Path Property.

As in case (a), this implies that after possibly subdividing at most two edges of $\Gamma$, the resulting graph $\Gamma'$ possesses a reduced path $p'$ such that $v = \mu(p')$ is a subword of $r_{ij}^*$ missing at most three letters and such that $p'$ is not a part of a relator cycle in $\Gamma'$. Thus $2m_{ij} - 3 \leq |v| \leq 2m_{ij}$.

**Claim 2.** The path $p'$ is a simple non-closed path in $\Gamma'$.

If not, then $p'$ contains a nontrivial subpath $p_1$ such that $p'$ is a simple closed path in $\Gamma'$. We have $|\mu(p_1)| \neq 1$ since by assumption $H$ does not contain conjugates of elements of $A$. Thus $|\mu(p_1)| \geq 2$. Since $\Gamma'$ is folded, the first and the last edges letters of $\mu(p_1)$ are distinct.

Without loss of generality we may assume that the first letter of $\mu(p_1)$ is $a_i$ and the last letter of $\mu(p_1)$ is $a_j$. Thus $u := \mu(p_1) = (a_i a_j)^z$, where $1 \leq z \leq m_{ij}$. Hence $p_1 = e_1 \ldots e_t$. Denote the first (which is also the last) vertex of $p_1$ by $x$. Since $\Gamma'$ is folded, the only edge incident to $x$ and with label starting with $a_i$ is the edge $e_1$. Similarly, the only edge incident to $x$ and with label starting with $a_j$ is $e_t^{-1}$. Since $v = \mu(p')$ is a subword of $r_{ij}^*$, this implies that $p'$ is contained in the subgraph $e_1 \cup \cdots \cup e_t$ of $\Gamma'$. If $z|m_{ij}$ then clearly $p'$ is a subpath in a relator-cycle in $\Gamma'$, contrary to our assumption. Thus $z$ does not divide $m_{ij}$. Let $d = gcd(z, m_{ij})$, so that $1 \leq d < z$. We attach to $\Gamma'$ a loop-edge $f$ based at $x$ labeled by $(a_i a_j)^d$. The resulting graph $\Gamma_1$ represents the same subgroup as $\Gamma$ but has Euler characteristic one less than that of $\Gamma$. We then perform all the folding moves necessary to "wrap" the simple circuit $p_1$ around $f$. The resulting graph $\Gamma_2$ is clearly represents a conjugate of $H$ and has Euler characteristic no smaller than that of $\Gamma$. Since $d < z$, the graph $\Gamma_1$ has smaller complexity than $\Gamma$, contradicting the choice of $\Gamma$.

Thus Claim 2 is established and $p'$ is indeed a simple non-closed path in $\Gamma'$. The rest of the argument proceeds exactly as in part (a) and we leave the details to the reader. $\qquad\square$

**Theorem A.** *Let $k \geq 2$ be an integer and let $G$ be a Coxeter group given by presentation* († ) *such that for all $i < j$ we have $m_{ij} \geq 3k + 1$. Then the following holds:*

1. *Every torsion-free $k$-generated subgroup of $G$ is free and quasiconvex in $G$.*
2. *Every $k$-generated subgroup of $G$, which does not contain any conjugates of the generators $a_i$, is quasiconvex in $G$.*
3. *Suppose that all finite $m_{ij}$ are even. Then every subgroup of $G$ generated by at most $(k+1)/2$ elements, is quasiconvex in $G$.*

*Proof.* The assumption that $m_{ij} \geq 3k + 1$ implies that $m_{ij} \geq 4$ and hence $G$ is word-hyperbolic. First observe that (2) implies (3). Indeed, suppose (2) is known to hold and $G$ is as in (3). Consider the map $\eta : A \to C = \langle c \,|\, c^2 = 1 \rangle$ defined as $\eta(a_i) = c$ for $i = 1, \ldots, n$. This map extends to a homomorphism $\eta : G \to C$ since all $m_{ij}$ are assumed to be even. Put $G_0 := ker(\eta)$. Then $G_0$ has index two in $G$ and $G_0$ does not contain any conjugate of $a_i$, $i = 1, \ldots, n$. Let $H \leq G$ be a subgroup generated by $s \leq (k+1)/2$ elements. Then $H_0 = H \cap G_0$ has index at most two in $H$ and hence by the Schreier formula [40] $H_0$ is generated by at most $2(s-1)+1 \leq k$ elements. By part (2) $H_0$ is quasiconvex in $G$. Since quasiconvexity is commensurability invariant (see for example [54, 43, 33]), this implies that $H$ is quasiconvex in $G$ as well.

We will now establish (1). Suppose now $H$ is a nontrivial $k$-generated and torsion-free subgroup of $G$. Since cyclic subgroups are free and always quasiconvex in hyperbolic groups [4, 14, 27], we may assume that $H$ is not cyclic. Since $H$ is $k$-generated, it can be represented by an $F$-graph with Euler characteristic $\geq 1-k$ (that is with fundamental group being a free group of rank at most $k$). Among all $F$-graphs with Euler characteristic $\geq 1-k$ representing subgroups conjugate to $H$ in $G$ choose a graph $\Gamma$ of minimal complexity. Therefore by Proposition 4.3 $\Gamma$ is folded, connected and every vertex in $\Gamma$ has degree at least three.

Choose a base-vertex $x_0$ in $\Gamma$. Recall that we have a homomorphism $\phi : \pi_1(\Gamma, x_0) \to G$ with the image of $\phi$ being a conjugate of $H$ in $G$. We claim that $\phi$ is a monomorphism (and hence $H$ is free). Suppose this is not the case. Then there exists a nontrivial reduced edge-path $p$ from $x_0$ to $x_0$ in $\Gamma$ such that $\overline{\mu(p)} = 1 \in \overline{w}$. Denote the label of $p$ by $w$. Since $\overline{w} = 1 \in G$, Proposition 3.1 implies that $w$ has a subword $v$ such that $v$ is also a subword of some $r_{ij}^*$ with $|r_{ij}| - 3 \leq |v|$. Thus $w$ is reduced but not weakly Dehn-reduced. However, this is impossible by Proposition 4.3 $\Gamma$. Thus $\phi$ is indeed a monomorphism and hence $H$ is free. Moreover, Proposition 4.3 implies that $\Gamma$ has the Relator Path Property. Hence by Proposition 4.2 $H$ is quasiconvex in $G$ (since quasiconvexity of a subgroup is preserved by conjugation).

The proof of part (2) is virtually identical to that of part (1). We choose $\Gamma$ exactly as in (1). Proposition 4.3 implies that $\Gamma$ has the Relator Path Property. Hence, again, by Proposition 4.2 $H$ is quasiconvex in $G$. This completes the proof of Theorem A.                                                    □

The following lemma is essentially Lemma 1.1 in Scott [53]:

**Lemma 4.4.** *Let $G$ be a finitely generated group and let $H \leq G$ be a finitely generated subgroup separable in $G$. Suppose $H \leq H_1 \leq G$, where $[H_1 : H] < \infty$. Then $H_1$ is separable in $G$.*

**Definition 4.5.** Let $G$ be a Coxeter group given by presentation ($\ddagger$) and suppose $H \leq G$ is a finitely generated subgroup and $g \in G - H$ is an element.

We say that a an $F$-graph $\Gamma$ is a *graph for the pair* $(H, g)$ if $\Gamma$ is a graph with two distinguished vertices $x_0, x_1$ such that:
(1) $\phi(\pi_1(\Gamma, x_0)) = H \leq G$
(2) For some path $p$ from $x_0$ to $x_1$ in $\Gamma$ we have $\overline{\mu(p)} = g \in G$.

Note that since $g \notin H$, we automatically have $x_1 \neq x_0$. It is also clear that a folding move on a graph for the pair $(H, g)$ produces a new graph for the pair $(H, g)$.

**Theorem B.** *Let $k \geq 2$ and $G$ be as in ($\dagger$) such that for each $i < j$ we have $m_{ij} \geq 3k+7$. Suppose in addition that ($\dagger$) satisfies the Separability Condition.*

*Then*

1. *Every $k$-generated subgroup of $G$, which does not contain any conjugate of the any generator $a_i$, is separable in $G$.*
2. *Every subgroup of $G$ generated by at most $(k+1)/2$ elements is separable in $G$.*

*Proof.* We again observe that (1) implies (2) as before. Indeed, suppose (1) is known to hold and $G$ is as in (2). Consider again the homomorphism $\eta : A \to C = \langle c \,|\, c^2 = 1 \rangle$ defined by $\eta(a_i) = c$ for $i = 1, \ldots, t$ and let $G_0 := ker(\eta)$. Then $G_0$ has index two in $G$ and $G_0$ does not contain any conjugate of $a_i$, $i = 1, \ldots, t$. If $H \leq G$ is a subgroup generated by $s \leq (k+1)/2$ elements then $H_0 = H \cap G_0$ has index at most two in $H$ and is hence generated by at most $2(s-1)+1 \leq k$ elements. Then by part (1) the subgroup $H_0$ is separable in $G$. Therefore by Lemma 4.4 $H$ is separable in $G$ as well, as required.

We will now establish (1). Let $H \leq G$ be a $k$-generated subgroup such that $H$ does not contain any conjugates of the generators $a_i$. Note that since $a_1 \notin H$, we have $G \neq H$. Let $g \in G - H$ be an arbitrary element. Note that if $\Gamma$ is an $F$-graph representing $H$ with base-vertex $x_0$, we can enlarge it to a graph for the pair $(H, g)$ as follows: take a non-loop edge labeled by a geodesic representative of $g$ and attach it to $\Gamma$ at the base-vertex $x_0$. Clearly the resulting graph is a graph for the pair $(H, g)$-pair having the same Euler characteristic as $\Gamma$.

Among all graphs for the pair $(H, g)$ with Euler characteristic $\geq 1 - k$ choose a graph $(\Gamma, x_0, x_1)$ with minimal complexity. Then $\Gamma$ is obviously folded, connected and has at most two vertices of degree one (namely $x_0, x_1$). Recall that $x_0 \neq x_1$ by definition.

We claim that $\Gamma$ has the Relator Path Property. Indeed, suppose not. Then after possibly subdividing two edges of $\Gamma$, the resulting graph $\Gamma'$ possesses a reduced path $p$ such that $v = \mu(p)$ is a subword of $r_{ij}^*$ missing at most three letters and such that $p$ is not a part of a relator cycle in $\Gamma'$. Thus $2m_{ij} - 3 \leq |v| \leq 2m_{ij}$. Note that

$\Gamma$ and $\Gamma'$ have the same Euler characteristics. Moreover, if $\sigma(\Gamma) = (l, q)$, $\sigma(\Gamma') = (l', q')$ then $l = l'$ and $q' \leq q + 2$.

**Claim.** The path $p$ is a simple path in $\Gamma'$.

If not, then $p$ contains a nontrivial subpath $p'$ such that $p'$ is a simple closed path in $\Gamma'$. We have $|\mu(p')| \neq 1$ since by assumption $H$ does not contain conjugates of elements of $A$. Thus $|\mu(p')| \geq 2$. Since $\Gamma'$ is folded, the first and the last letters of $\mu(p')$ are distinct. Without loss of generality we may assume that the first letter of $\mu(p')$ is $a_i$ and the last letter is $a_j$. Thus $u := \mu(p') = (a_i a_j)^z$, where $1 \leq z \leq m_{ij}$. Let $p'$ be of the form $p' = e_1 \ldots e_t$ where the $e_h$ are edges of $\Gamma'$. Denote the first (which is also the last) vertex of $p'$ by $x$. Since $\Gamma'$ is folded, the only edge incident to $x$ with label starting with $a_i$ is $e_1$. Similarly, the only edge incident to $x$ with label starting with $a_j$ is $e_t^{-1}$. Since $v = \mu(p)$ is a subword of $r_{ij}^*$, this implies that $p$ is contained in the subgraph $e_1 \cup \cdots \cup e_q$ of $\Gamma'$. It $z | m_{ij}$ then clearly $p$ is a subpath of a relator-cycle in $\Gamma'$, contrary to our assumption. Thus $z$ does not divide $m_{ij}$. Let $d := gcd(z, m_{ij})$, so that $1 \leq d < z$. We attach to $\Gamma$ a loop-edge $f$ at $x$ labeled by $(a_i a_j)^d$. The resulting graph $\Gamma_1$ is still an $(H, g)$-pair of Euler characteristic one less than that of $\Gamma$. We then perform all the folding moves necessary to "wrap" the simple circuit $p'$ around $f$. The resulting graph $\Gamma_2$ is clearly an $(H, g)$-pair of Euler characteristic no smaller than that of $\Gamma$. Since $d < z$, the graph $\Gamma_2$ has smaller complexity than $\Gamma$, contradicting the choice of $\Gamma$.

Thus the Claim is established and $p$ is indeed a simple non-closed path in $\Gamma'$. Recall that $\Gamma$ has at most two vertices of degree one and $\Gamma'$ is obtained from $\Gamma$ by subdividing at most two edges. Thus by possibly removing the "spikes" corresponding to degree-one vertices and then performing at most four inverse subdivisions on $\Gamma'$, we obtain a graph of the same Euler characteristic as $\Gamma$ and with every vertex of degree at least three. Therefore Lemma 2.13 implies that $p$ consists of at most $(2k - 3) + 6 = 2k + 3$ edges, provided $\pi_1(\Gamma, x_0)$ is non-cyclic. If $\pi_1(\Gamma, x_0)$ is cyclic, then $\Gamma$ has at most four edges. In this case it is easy to see that $p$ consists of at most five edges. Since $k \geq 2$ and $5 < 2k + 3$, we conclude that in either case $|p| \leq 2k + 3$. Therefore there is an edge $e$ of $p$ such that $|\mu(e)| \geq |\mu(p)|/(2k + 3)$.

Recall that $|\mu(p)| = |v| \geq 2m_{ij} - 3$. Hence

$$|\mu(e)| \geq |v|/(2k + 3) \geq \frac{2m_{ij} - 3}{2k + 3} > 3,$$

where the last inequality holds since by assumption $m_{ij} \geq 3k + 7$. Denote the first vertex of $p$ by $x$ and the last vertex of $p$ by $y$.

We now perform a move of type $A0$ by adding an edge from $x$ to $y$ with label of length at most three corresponding to the part of $r_{ij}^*$ missing in $v$. Next we perform a type-$A1$ move and remove the edge $e$ of length at least four. Note that the vertex $x_1$ belongs to the resulting connected graph $\Gamma_1$. Therefore $\Gamma_1$ is an $(H, g)$-pair of the same Euler characteristic as $\Gamma$ but with smaller complexity. This contradicts the minimal choice of $\Gamma$.

Thus $\Gamma$ indeed has the Relator Path Property, as claimed. By Proposition 4.2 there exists a homomorphism $\psi : G \to K$, where $K$ is a finite group (which can be taken to be the symmetric group on the set of vertices of the 2-completion) such that $\psi(g) \notin \psi(H)$. Since $g \in G - H$ was chosen arbitrarily, this implies that $H$ is separable in $G$, as required. □

## 5. Small cancellation theory and Artin groups

Our next goal is to prove

**Theorem C.** *Let $k \geq 2$ be an integer. Let*

($\spadesuit$) $$G = \langle a_1, \ldots, a_n \,|\, u_{ij} = u_{ji}, \; where \; 1 \leq i < j \leq t \rangle$$

be an Artin group where $|u_{ij}| = m_{ij} \geq 5k$ for each $i < j$. Suppose $H \leq G$ is a $k$-generated subgroup such that $H$ has trivial intersection with any conjugate of every two-generator subgroup $G_{ij} = \langle a_i, a_j \rangle$ for which $m_{ij} < \infty$.

Then $H$ is free.

We need to recall some facts from Appel and Schupp [7]. about how to apply small cancellation theory to Artin groups. They first considered the two-generator one-relator Artin group

$$G_{ij} = \langle a_i, a_j \,|\, u_{ij} = u_{ji} \rangle,$$

and showed that the symmetrized set generated by the defining relator satisfies the small cancellation conditions $C(4)$ and $T(4)$. They also introduced *strips* to study the finer geometry of $C(4) - T(4)$ diagrams. ( "Strips" have also been called "chains" or "corridors" in some later articles.) Now for the free group $F(A) = F(a_1, \ldots, a_n)$ every nontrivial reduced word $w$ has a unique *normal form with exponents*

$$w = a_{h_1}^{n_1} \ldots a_{h_s}^{n_s}$$

where each $j_t \neq j_{t+1}$ and each $n_h \neq 0$. The integer $s$ is the *syllable length* of $w$ and is written $||w||$. (This is the way one would write elements if we considered $F$ as the free product of the infinite cyclic groups on the generators, but we are not using the free product structure and we just consider elements of the free group as being written in normal form with exponents.)

For each $i < j$ such that $m_{ij} < \infty$ denote by $R_{ij}$ the set of all freely reduced and cyclically reduced words in $F(a_i, a_j)$ which are equal to one in the group $G_{ij}$. A consequence of the theory is that the canonical homomorphism from $G_{ij}$ into $G$ is injective [7] provided all $m_{ij} \geq 4$ in ($\ddagger$).

The results about strips [7] prove:

**Proposition 5.1.** *Suppose $2 \leq m_{ij} < \infty$ and let $w$ be a nonempty word from $R_{ij}$. Then $||w|| \geq 2m_{ij}$.*

We will work with the following infinite presentation of the Artin group $G$:

($\clubsuit$) $$G = \langle a_1, \ldots, a_n \,|\, \mathcal{R} \rangle.$$

where

$$\mathcal{R} = \bigcup_{\substack{i < j \\ m_{ij} < \infty}} R_{ij}.$$

For the remainder of this section we will assume that $k \geq 2$ and that $G$ is an Artin group satisfying the assumptions of Theorem C. The notations described above and presentation ($\clubsuit$) will also be fixed for the remainder of this section.

The point of shifting to the infinite presentation ($\clubsuit$) is that it allows a strong use of minimality. In considering van Kampen diagrams showing that a word $w$ is a consequence of the defining relators one need only consider diagrams which are minimal for $w$, that is, diagrams with as few regions as possible over all $\mathcal{R}$-diagrams for $w$. In such a diagram, distinct regions labeled by relators from the same set $R_{ij}$ cannot have even a vertex in common for otherwise they could be combined into a single region, contradicting minimality. So if two regions in a minimal diagram have an edge in common, the label on that edge can only be a power of a single generator, say $a_i$, and thus, in view of the proposition above and the fact that all $m_{ij} \geq 4$, an interior region in a minimal diagram must have degree at least eight. The geometry of minimal diagrams for the infinite presentation is thus essentially that of a set of relators satisfying $C(8)$ and the following proposition holds [7].

**Proposition 5.2.** *If $w$ is a nontrivial freely reduced word representing $1$ in $G$ then $w$ contains a subword $v$ such that $v$ is also a subword of some $r \in R_{ij}$ with $r = vu$, $||u|| \leq 3$ and $||v|| \geq 2m_{ij} - 3$.*

In order to study subgroups of Artin groups we need to incorporate syllable length into our measure of the complexity of subgroup graphs. The difference between the previous definition and the definition below is that syllable length becomes the most important component.

**Definition 5.3** (Fine Complexity). Let $\Gamma$ be an $F$-graph. As before, we partition the set $E(\Gamma)$ of oriented edges of $\Gamma$ as $E(\Gamma) = E^+ \cup E^-$ where for any $e \in E(\Gamma)$ we have $e \in E^+ \iff e^{-1} \in E^-$ and where $E^+ \cap E^- = \emptyset$.

Let $s = \sum_{e \in E^+} ||\mu(e)||$, let $l = \sum_{e \in E^+} |\mu(e)|$ and let $q$ be the number of vertices in $\Gamma$. We now define the *fine complexity* $c(\Gamma)$ to be $c(\Gamma) = (s, l, q)$.

We order complexities lexicographically as usual, that is, $(s, l, q) < (s', l', q')$ if one of the following occurs:

1. $s < s'$, or
2. $s = s'$ and $l < l'$; or
3. $s = s'$ and $l = l'$ and $q < q'$.

**Lemma 5.4.** *Let $H \leq G$ be a $k$-generated non-cyclic subgroup. (Hence $k \geq 2$ and $H$ can be represented by an $F$-graph with Euler characteristic $\geq 1 - k$.)*

*Among all connected $F$-graphs with Euler characteristic $\geq 1 - k$ representing conjugates of $H$ in $G$ choose a graph $\Gamma$ of minimal fine complexity.*

*Then $\Gamma$ is folded and the degree of every vertex in $\Gamma$ is at least three.*

*Proof.* It is clear that if $\Gamma$ has a degree-one vertex, then there is a graph representing a conjugate of $H$ and with the same Euler characteristic as $\Gamma$ but with smaller fine complexity. Also, suppose $x$ is a vertex of $\Gamma$ of degree two. If there is a loop-edge at $x$ then $\Gamma$ consists of a single loop-edge, contradicting the assumption that $H$ is not cyclic. So there are two distinct non-loop edges $e_1, e_2$ originating at $x$ with $\mu(e_i) = w_i$ and such that there are no other edges incident to $x$. If $w_1 = w_2^{-1}$, we can fold these two edges into a "spike" and obtain a graph of obviously smaller fine complexity yielding a contradiction. Suppose now that $w_1 \neq w_2$. Let $w$ be the reduced form of $w_1^{-1}w_2$. We can perform move $A3$, namely remove the edges $e_2, e_1$ and the vertex $x$ from $\Gamma$ and insert a new edge $e$ with label $w$ going from the terminal vertex of $e_1$ to the terminal vertex $e_2$. Obviously, the resulting graph $\Gamma'$ represents a conjugate of $H$. Also the Euler characteristics of $\Gamma$ and $\Gamma'$ are the same. Note that $||w_1^{-1}w_2|| \leq ||w_1|| + ||w_2||$ and $|w_1^{-1}w_2| \leq |w_1| + |w_2|$. However $\Gamma'$ has fewer vertices than $\Gamma$ and hence $\Gamma'$ has lower fine complexity that $\Gamma$, contrary to our assumptions.

Thus we have shown that the degree of every vertex in $\Gamma$ is at least three.

Suppose now that $\Gamma$ is not folded. Thus there are edges $e_1 \neq e_2$ of $\Gamma$ with labels $w_1$ and $w_2$ accordingly and with a common initial vertex $x$ such that $w_1^{-1}w_2$ is not freely reduced.

Let $v$ be the maximal common initial segment of $w_1$ and $w_2$.

Thus we have graphic equalities $w_1 = vu_1$ and $w_2 = vu_2$ where the product $u_1^{-1} \cdot u_1$ is reduced. Let $\Gamma'$ be obtained from $\Gamma$ by a fold corresponding to this situation. Let $c(\Gamma) = (s, l, q)$ and $c(\Gamma') = (s', l', q')$. By the choice of $v$ for at least one of $i = 1, 2$ we have $||w_i|| = ||v|| + ||u_i||$ and for each $i = 1, 2$ $||v|| + ||u_i|| - 1 \leq ||w_i|| \leq ||v|| + ||u_i||$. Also $||v|| \geq 1$. Hence

$$||w_1|| + ||w_2|| \geq 2||v|| + ||u_1|| + ||u_2|| - 1 \geq ||v|| + ||u_1|| + ||u_2||,$$

which implies that $s \geq s'$.

Clearly $l' < l$ since $|v| > 0$. Thus $c(\Gamma') < c(\Gamma)$ which contradicts the choice of $\Gamma$. Thus $\Gamma$ is folded and the lemma is proved. $\qquad\square$

*Proof of Theorem C.* Suppose $H$ is a nontrivial $k$-generated subgroup of $G$. If $H$ is cyclic then it is infinite cyclic since $G$ is torsion-free by the result of [7]. The statement of Theorem C clearly holds in this case.

Thus we assume that $H$ is not cyclic. Among all connected $F$-graphs with Euler characteristic $\geq 1 - k$ representing conjugates of $H$ in $G$ choose a graph $\Gamma$ with minimal fine complexity. By Lemma 5.4 $\Gamma$ is folded and every vertex in $\Gamma$ has degree at least three. Hence the label of every reduced edge-path in $\Gamma$ is a freely reduced word. Let $c(\Gamma) = (s, l, q)$.

Let $x_0$ be a base-vertex of $\Gamma$. Then we have a homomorphism $\phi : \pi_1(\Gamma, x_0) \to G$ with $\phi(\pi_1(\Gamma, x_0))$ being a conjugate of $H$ in $G$. We claim that $\phi$ is a monomorphism (and hence $H$ is free).

Suppose not. Then there exists a nontrivial reduced edge-path $p$ from $x_0$ to $x_0$ in $\Gamma$ such that $\overline{\mu(p)} = 1 \in G$. Denote the label of $p$ by $w$. Since $\overline{w} = 1 \in G$, Proposition 5.2 implies that $w$ has a subword $v$ such that $v$ is also a subword of some $r \in R_{ij}$ with $r = vu$, where $\|u\| \leq 3$, $\|v\| \geq \|r\| - 3 \geq 2m_{ij} - 3$.

We now inspect how the occurrence of $v$ in $w$ corresponds to traversing the path $p$. There is a subpath $p_1 = e_1, \ldots, e_t$ of $p$ with the following properties:

$$\mu(e_1) = y'y_1, \mu(e_t) = y_ty'', \mu(e_h) = y_h, 1 < h < t,$$

and there is a graphic equality

$$v = y_1 y_2 \ldots y_t,$$

and $|y_1| > 0$, $|y_t| > 0$. Note that we allow for the possibility when one or both $y'$, $y''$ are trivial.

For those of $y'$, $y''$ which are not equal to 1, we perform an $A2$-subdivision of the edges $e_1$, $e_t$ accordingly into two edges with labels $y'$, $y_1$ and $y_t$, $y''$ accordingly. Denote the resulting graph by $\Gamma'$ and let $c(\Gamma') = (s', l', q')$. Note that $s' \leq s + 2$, $l' = l$ and $q' \leq q + 2$.

Let $e_1'$ be the edge with label $y_1$ and $e_t'$ be the edge with label $y_t$ resulting from this operation. Denote $p_1' = e_1', e_2, \ldots, e_{t-1}, e_t'$.

**Claim.** The path $p_1' = e_1', e_2, \ldots, e_{t-1}, e_t'$ is a simple path in $\Gamma'$.

Indeed, suppose this is not the case. Recall that $p_1'$ is a reduced path. Choose a nonempty simple closed subpath $p_1''$ of $p_1'$. Notice that $p_1''$ is in fact a simple closed path in the original graph $\Gamma$ unless $p_1'' = p_1'$.

Denote $u = \mu(p_1'')$. Suppose first that $\overline{u} =_G 1$. If $p_1'' \neq p_1'$, then we can perform move $A1$ on the original graph $\Gamma$ and remove the first edge of $p_1''$. The resulting graph represents the same subgroup as $\Gamma$ but has lower fine complexity, contradicting our assumptions. Suppose now $p_1'' = p_1'$, so that $p_1'$ is a simple closed path in $\Gamma'$. If $\Gamma = \Gamma'$ then we can perform an $A1$ move on $\Gamma$ and obtain a contradiction as before. Thus $\Gamma \neq \Gamma'$, which means that $e_1 = e_t = e$ and $\mu(e) = y_ty_1$, so that $y' = y_t$ and $y'' = y_1$ are both nonempty. Then $s \leq s' \leq s + 1$. If $s' = s$ (that is $\|y_ty_1\| = \|y_t\| + \|y_1\|$), then we can perform an $A1$-move on $\Gamma'$ deleting the edge $e_1'$ with label $y_1$ and obtain a graph of lower fine complexity, yielding a contradiction. Suppose now that $s' = s + 1$. If $p_1'$ consists of at least three edges, then we can still perform the move $A1$ on $\Gamma'$ and remove the second edge of $p_1'$. The resulting graph has smaller fine complexity than $\Gamma$, again giving us a contradiction. Suppose now that $t = 2$ and $p_1' = e_1' e_2'$. In this case $e = e_1 = e_2$ is a loop-edge of $\Gamma$ with label $\mu(e) = y_2 y_1$. We know that $u = y_1 y_2 =_G 1$ and therefore $\mu(e) =_G 1$. Again, we can perform move $A1$ on $\Gamma$ and remove $e$, thus lowering the fine complexity and obtaining a contradiction.

Thus $\overline{u} \neq_G 1$. However, since $u$ is a subword of $R \in R_{ij}$ and $u$ is the label of a circuit in $\Gamma'$ , this implies that a conjugate of $H$ contains a nontrivial element of $G_{ij} = \langle a_i, a_j \rangle$, contrary to the assumptions of Theorem C. Thus the Claim is verified.

We now know that $p_1'$ is a simple non-closed path in $\Gamma'$. Lemma 2.13 and Lemma 5.4 imply that $p_1'$ has at most $(2k - 3) + 2 = 2k - 1$ edges and so $t \leq 2k - 1$. Recall that $v = y_1 y_2 \ldots y_t$.

Therefore for some $1 \leq h \leq t$ we have $\|y_h\| \geq \|v\|/(2k-1)$. Recall also that $\|v\| \geq \|r\| - 3 \geq 2m_{ij} - 3$. Hence

$$\|y_h\| \geq \|v\|/(2k-1) = (2m_{ij} - 3)/(2k - 1) > 5,$$

where the last inequality holds since by assumption $m_{ij} \geq 5k$.

We now perform a move of type $A0$ on $\Gamma'$ by adding an edge labeled $u$ from the terminal vertex of $e_t'$ to the initial vertex of $e_1'$. Recall that $R = vu$ and $\|u\| \leq 3$. Next we perform a type-$A1$ move and remove the $h$-th edge of $p_1'$ (labeled $y_h$). By Lemma 5.4 the resulting graph $\Gamma_1$ represents a subgroup conjugate to $H$ in $G$. It is also clear that $\Gamma_1$ has the same Euler characteristic as $\Gamma$ and $\Gamma'$.

Since $s' \leq s + 2$, $||u|| \leq 3$ and $||y_h|| > 5$, the graph $\Gamma_1$ has smaller fine complexity than does $\Gamma$, which contradicts the choice of $\Gamma$. Thus we have shown that $\phi$ is a monomorphism and hence $H$ is free. This completes the proof of Theorem C. $\qquad\blacksquare$

## 6. FREE SUBGROUPS OF ONE-RELATOR GROUPS WITH TORSION

We recall the following strengthened version of Newman's "Spelling Theorem" [45, 40].

**Theorem 6.1.** *Let $G$ be given by a presentation*

(*) $$G = \langle a_1, \ldots, a_n \,|\, r^m = 1 \rangle,$$

*where $m > 1$ and where $r$ is a nontrivial cyclically reduced word in $F = F(a_1, \ldots, a_n)$ which is not a proper power. Suppose $w$ is a freely reduced word in $F$ such that $w =_G 1$. Then $w$ contains a subword $v$ of the form $v = r_*^{m-1} x$ where $r_*$ is a cyclic permutation of $r$ or $r^{-1}$ and where $x$ is an initial segment of $r_*$ containing an occurrence of every letter $a_i$ which appears in $r$.*

This strengthened version of Newman's Spelling Theorem was proven by Gurevich [29] and used by Pride in his proof that all two-generator torsion-free subgroups of one-relator groups with torsion are free. A geometric proof using diagrams was given by Schupp[51] who also strengthened the statement of the Frieheitssatz in the case where the relator is not a proper power.

We shall also need the well-known description of elements of finite order in one-relator groups with torsion, which is due to Karrass, Magnus and Solitar [36]:

**Proposition 6.2.** *Let $G$ be as in Theorem 6.1. Then the word $r$ represents an element of order $m$ in $G$. Moreover, any nontrivial element of finite order in $G$ is conjugate to $r^d$ where $0 < |d| < m$.*

**Convention 6.3.** For the remainder of the article we set $F = F(A) = F(a_1, \ldots, a_n)$ and assume that $G$ is a group given by presentation (*). We return to the same definition of the complexity of a subgroup used in the section on Coxeter groups.

The following proposition is essentially Lemma 4 of Lyndon-Schützenberger [41].

**Proposition 6.4.** *Let $r$ be a nontrivial reduced and cyclically reduced word which is not a proper power in $F$. Let $u$ be a reduced word with $|u| \geq |r|$. Suppose $r^d = xuy = x'uy'$ for some $d > 0$ and $x' = xw$ (where all equalities are graphical). Then $w = r_*^z$ where $0 \leq z \leq d$ and $r_*$ is a cyclic permutation of $r$.*

*Moreover, $r^{-1}$ is not a subword of $r^d$ for any $d > 0$.*

**Definition 6.5.** Suppose $r$ is a cyclically reduced word which is not a proper power in $F(A)$. For each $a \in A \cup A^{-1} = \{a_1, \ldots, a_n\}^{\pm 1}$ we denote by $b_a(r)$ the total number of times the letter $a$ occur in $r$ (not counting the occurrences of $a^{-1}$). The *letter bound*, $b(r)$ of $r$ is

$$b(r) := \max\{b_a(r) \,|\, a \in A \cup A^{-1}\}.$$

**Theorem D.** *Let $k \geq 2$ be an integer and let*

(3) $$G = \langle a_1, \ldots, a_n \,|\, r^m = 1 \rangle$$

*be a one-relator group where $r$ is a nontrivial cyclically reduced word which is not a proper power and let $b = b(r)$. If $m \geq b(6k - 4) + 2$ then every $k$-generated torsion-free subgroup of $G$ is free.*

*Proof.* If $|r| = 1$, then $G$ is virtually free and Theorem D obviously holds. Thus we may assume that $|r| \geq 2$.

Suppose $H$ is a nontrivial torsion-free $k$-generated subgroup of $G$. If $H$ is cyclic then it obviously free and there is nothing to prove. Thus we assume that $H$ is not cyclic. Among all connected $F$-graphs with Euler characteristic $\geq 1 - k$ representing conjugates of $H$ in $G$ choose a graph $\Gamma$ with minimal complexity. It is clear that $\Gamma$ is folded and has no degree-one vertices. Hence the label of every reduced edge-path in $\Gamma$ is a freely reduced word. Moreover, since $H$ is not cyclic and because of the minimality assumption, $\Gamma$ has no degree-two vertices.

Let $x_0$ be a base-vertex of $\Gamma$. Then we have a homomorphism $\phi : \pi_1(\Gamma, x_0) \to G$ with $\phi(\pi_1(\Gamma, x_0))$ being a conjugate of $H$ in $G$. We claim that $\phi$ is a monomorphism (and hence $H$ is free).

Suppose not. Then after possibly subdividing at most two edges of $\Gamma$ the resulting graph $\Gamma'$ possesses a nontrivial reduced edge-path $p$ such that $v = \overline{\mu(p)}$ has the form $v = r_*^{m-1}a$ where $r_*$ is a cyclic permutation of $r$ and where $a$ is the first letter of $r_*$.

Let

$$p = e_1 e_2 \ldots e_{t-1} e_t$$

where each $e_i$ is an edge of $\Gamma'$. Denote the initial vertex of $p$ by $x$ and the terminal vertex of $p$ by $y$.

**Claim 1.** Suppose some $e_i$ is disjoint from all the edges $e_j, j \neq i$. Then $|\mu(e_i)| < |r|$.

Indeed, suppose $|\mu(e_i)| \geq |r|$. Recall that $r_* = a\beta$ where $|\beta| = |r| - 1$ and $r_*^n = v\beta$. Recall also that $x$ and $y$ are the initial and the terminal vertices of $p$. We first perform move $A0$ on $\Gamma'$ and attach an edge labeled $\beta$ going from $y$ to $x$. We then perform move $A1$ and remove the edge $e_i$. The resulting graph $\Gamma''$ represents the same subgroup as $\Gamma$ and the Euler characteristics of $\Gamma$ and $\Gamma''$ are equal. However $\Gamma''$ has smaller complexity than $\Gamma$ yielding a contradiction.

**Claim 2.** Suppose that for some $1 \leq i < j \leq s$ we have $e_i = e_j^{\pm 1}$. Then $|\mu(e_i)| < |r|$.

Suppose, on the contrary, that $|\mu(e_i)| \geq |r|$. By Proposition 6.4 this implies that $e_i = e_j$. Moreover, since $|\mu(e_i)| \geq |r|$, by Proposition 6.4 the label of the path $e_i \ldots e_{j-1}$ has the form $(r')^z$ where $r'$ is some cyclic permutation of $r$ and $1 \leq z < m$. Since the path $e_i \ldots e_{j-1}$ is a cycle in $\Gamma'$, this implies that $H$ contains an element conjugate to $(r')^z$, which contradicts the assumption that $H$ is torsion-free. Thus Claim 2 is established.

We now return to the proof of Theorem D. Claim 1 and Claim 2 imply that $|\mu(e_i)| < |r|$ for $i = 1, \ldots, s$. Recall that $\mu(p) = v = r_*^{m-1}a$ where $r_*$ is a cyclic permutation of $r$ and where $a$ is the first letter of $r_*$. Therefore $s > n - 1$, that is $s \geq n$.

By Lemma 2.13 the maximal tree in $\Gamma'$ has at most $(2k-3)+2 = 2k-1$ non-oriented edges. Since the rank of the fundamental group of $\Gamma'$ is bounded above by $k$, the graph $\Gamma'$ has at most $2k-1+k = 3k-1$ non-oriented edges and at most $6k-2$ oriented edges.

Recall that $t > m > b(r)(6k-2)+1$. Therefore some oriented edge $e$ of $\Gamma'$ is repeated in $p = e_1 \ldots e_{t-1}e_t$ at least $b(r)+1$ times. Therefore by definition of $b(r)$ there are $i < j$ such that $e_i = e = e_j$ and the label of $e_i \ldots e_{j-1}$ has the form $(r')^z$ where $0 < z < m$ and $r'$ is a cyclic permutation of $r$. Since $e_i \ldots e_{j-1}$ is a cycle in $\Gamma'$, this implies that $H$ contains an element conjugate to $r^z$, contradicting our assumption that $H$ is torsion-free. This completes the proof of Theorem D. $\qquad\square$

## References

[1] G. Arzhantseva, *On groups in which subgroups with a fixed number of generators are free,* (Russian) Fundam. Prikl. Mat. **3** (1997), no. 3, 675–683

[2] G. Arzhantseva, *Generic properties of finitely presented groups and Howson's theorem,* Comm. Algebra **26** (1998), no. 11, 3783–3792

[3] G. Arzhantseva, *A property of subgroups of infinite index in a free group,* Proc. Amer. Math. Soc. **128** (2000), no. 11, 3205–3210

[4] J.Alonso, T.Brady, D.Cooper, V.Ferlini, M.Lustig, M.Mihalik, M.Shapiro and H.Short, *Notes on hyperbolic groups,* In: " Group theory from a geometrical viewpoint", Proceedings of the workshop held in Trieste, É. Ghys, A. Haefliger and A. Verjovsky (editors). World Scientific Publishing Co., 1991

[5] I. Agol, D. Long and A. Reid, *The Bianchi groups are separable on geometrically finite subgroups,* Ann. of Math. (2) **153** (2001), no. 3, 599–621

[6] G. Arzhantseva and A. Olshanskii, *Genericity of the class of groups in which subgroups with a lesser number of generators are free,* (Russian) Mat. Zametki **59** (1996), no. 4, 489–496

[7] K. Appel and P. Schupp, *Artin groups and infinite Coxeter groups,* Invent. Math. **72** (1983), no. 2, 201–220

[8] M. Bestvina, *Non-positively curved aspects of Artin groups of finite type,* Geom. Topol. **3** (1999), 269–302

[9] I. Bumagin, *On Small Cancellation k-Generated Groups with (k − 1)-Generated Subgroups All Free,* Intern J. Alg. Comp. **11** (2001), no. 5, 507–524

[10] M. Bestvina and N. Brady, *Morse theory and finiteness properties of groups,* Invent. Math. **129** (1997), no. 3, 445–470

[11] K.-U. Bux and C. Gonzalez, *The Bestvina-Brady construction revisited: geometric computation of Σ-invariants for right-angled Artin groups,* J. London Math. Soc. (2) **60** (1999), no. 3, 793–801

[12] J. Becker, M. Horak and L. VanWyk, *Presentations of subgroups of Artin groups,* Missouri J. Math. Sci. **10** (1998), no. 1, 3–14

[13] J. Crisp, *Symmetrical subgroups of Artin groups,* Adv. Math. **152** (2000), no. 1, 159–177

[14] M. Coornaert, T. Delzant, and A. Papadopoulos, *Géométrie et théorie des groupes. Les groupes hyperboliques de Gromov.* Lecture Notes in Mathematics, 1441; Springer-Verlag, Berlin, 1990

[15] D. Cooper, D. Long and A. Reid, *Infinite Coxeter groups are virtually indicable,* Proc. Edinburgh Math. Soc. (2) **41** (1998), no. 2, 303–313

[16] R. Charney and M. Davis, *When is a Coxeter system determined by its Coxeter group?* J. London Math. Soc. (2) **61** (2000), no. 2, 441–461

[17] J. Crisp and L. Paris, *The solution to a conjecture of Tits on the subgroup generated by the squares of the generators of an Artin group,* Invent. Math. **145** (2001), no. 1, 19–36

[18] M. Davis and T. Januszkiewicz, *Right-angled Artin groups are commensurable with right-angled Coxeter groups,* J. Pure Appl. Algebra **153** (2000), no. 3, 229–235

[19] W. Dicks and I. Leary, *On subgroups of Coxeter groups.* in "Geometry and cohomology in group theory (Durham, 1994)", 124–160, London Math. Soc. Lecture Note Ser., **252**, Cambridge Univ. Press, Cambridge, 1998

[20] W. Dicks and I. Leary, *Presentations for subgroups of Artin groups,* Proc. Amer. Math. Soc. **127** (1999), no. 2, 343–348

[21] D. Epstein, J. Cannon, D. Holt, S. Levy, M. Paterson, W. Thurston, *Word Processing in Groups,* Jones and Bartlett, Boston, 1992

[22] R. Gitik, *On the profinite topology on negatively curved groups,* J. Algebra **219** (1999), no. 1, 80–86

[23] R. Gitik, *Doubles of groups and hyperbolic LERF 3-manifolds,* Ann. of Math. (2) **150** (1999), no. 3, 775–806

[24] R. Gitik, *Tameness and geodesic cores of subgroups,* J. Austral. Math. Soc. Ser. A **69** (2000), no. 2, 153–16

[25] R. Gitik, *On the profinite topology on Coxeter groups,* IHES preprint, 2001

[26] M. Gromov, *Hyperbolic Groups,* in "Essays in Group Theory (G.M.Gersten, editor)", MSRI publ. **8**, 1987, 75–263

[27] E. Ghys and P. de la Harpe (editors), *Sur les groupes hyperboliques d'aprés Mikhael Gromov,* Birkhäuser, Progress in Mathematics series, vol. **83**, 1990.

[28] S. Gersten and H. Short, *Rational subgroups of biautomatic groups,* Ann. Math. (2) **134** (1991), no. 1, 125–158

[29] G. A. Gurevich *On the conjugacy problem for one-relator groups,* Dokladi Akademii Nauk SSSR, Ser. Mat. **207** (1972), 1436–1439.

[30] C. Hruska and D. Wise *Towers, ladders, and the B.B.Newman spelling theorem,* to appear in Jour. Aust. Math. Soc.

[31] T. Hsu and D. Wise, *Separating quasiconvex subgroups of right-angled Artin groups,* to appear in Math. Z.

[32] I. Kapovich and A. Myasnikov, *Stallings foldings and the subgroup structure of free groups,* to appear in the Journal of Algebra

[33] I. Kapovich, and H. Short, *Greenberg's theorem for quasiconvex subgroups of word hyperbolic groups,* Canad. J. Math. **48** (1996), no. 6, 1224–1244

[34] I. Kapovich and Richard Weidmann, *Nielsen methods and groups acting on hyperbolic spaces,* preprint, 2000

[35] I. Kapovich and Richard Weidmann, *Freely indecomposable groups acting on hyperbolic spaces,* preprint, 2001

[36] A. Karrass, W. Magnus, and D. Solitar, *Elements of finite order in groups with a single defining relation,* Comm. Pure Appl. Math. **13** (1960), 57–66

[37] D. Long and G. Niblo *Subgroup separability and 3-manifold groups,* Math. Z. **207** (1991), no. 2, 209–215

[38] D. Long and A. Reid, *The fundamental group of the double of the figure-eight knot exterior is GFERF,* Bull. London Math. Soc. **33** (2001), no. 4, 391–396

[39] D. Long and A. Reid, *On Subgroup Separability in Hyperbolic Coxeter Groups,* Geometriae Dedicata **87** (2001), 245–260

[40] R. Lyndon and P. Schupp, *Combinatorial Group Theory,* Springer-Verlag, 1977. Reprinted in the Classics in Mathematics series, 2000.

[41] R. Lyndon and M. Schützenberger, *The equation $a^M = b^N c^P$ in a free group,* Michigan Math. J. **9** (1962), 289–298

[42] J. Meier, H. Meinert and L. VanWyk, *On the Σ-invariants of Artin groups,* Topology Appl. **110** (2001), no. 1, 71–81

[43] M. Mihalik and W. Towle, *Quasiconvex subgroups of negatively curved groups,* J. Pure Appl. Algebra **95** (1994), no. 3, 297–301

[44] J. McCammond and D. Wise, *Coherence, local quasiconvexity and the perimeter of 2-complexes,* preprint

[45] B. B. Newman, *Some results on one-relator groups,* Bull. Amer. Math. Soc. **74** (1968), 568–571

[46] G. Niblo and D. Wise, *Subgroup separability, knot groups and graph manifolds,* Proc. Amer. Math. Soc. **129** (2001), no. 3, 685–693

[47] L. Paris, *Parabolic subgroups of Artin groups,* J. Algebra **196** (1997), no. 2, 369–399

[48] S. Pride, *The two-generator subgroups of one-relator groups with torsion,* Trans. Amer. Math. Soc. *234* (1977), no. 2, 483–496

[49] S. Pride, *Small cancellation conditions satisfied by one-relator groups,* Math. Z. **184** (1983), no. 2, 283–286

[50] L. Reeves, *Rational subgroups of cubed 3-manifold groups,* Michigan Math. J. **42** (1995), no. 1, 109–126

[51] P. Schupp, *A strengthened Freiheitssatz,* Math. Ann.**221** (1976), no. 1, 73–80

[52] P. Schupp, *Coxeter groups, perimeter reduction, 2-completion and subgroup separability,* to appear in Geom. Dedicata

[53] P. Scott, *Subgroups of surface groups are almost geometric,* J. London Math. Soc.(2) **17** (1978), 555–565

[54] H. Short, *Quasiconvexity and a theorem of Howson's,* in "Group theory from a geometrical viewpoint (Trieste, 1990)", 168–176, World Sci. Publishing, River Edge, NJ, 1991

[55] J.R. Stallings, *Topology of finite graphs,* Invent. Math. **71** (1983), no. 3, 551–565

[56] R. Strebel, *Small cancellation groups,* Appendix in "Sur les groupes hyperboliques d'aprés Mikhael Gromov", E. Ghys and P. de la Harpe (editors), Birkhäuser, Progress in Mathematics series, vol. **83**, 1990.

[57] G. A. Swarup, *Geometric finiteness and rationality,* J. Pure Appl. Algebra **86** (1993), no. 3, 327–333

[58] G. A. Swarup, *Proof of a weak hyperbolization theorem,* Q. J. Math. **51** (2000), no. 4, 529–533

[59] D. Wise, *Subgroup separability of graphs of free groups with cyclic edge groups,* Q. J. Math. **51** (2000), no. 1, 107–129

[60] D. Wise, *Subgroup separability of the figure 8 knot group,* preprint

DEPARTMENT OF MATHEMATICS, UNIVERSITY OF ILLINOIS AT URBANA-CHAMPAIGN, 1409 WEST GREEN STREET, URBANA, IL 61801, USA

*E-mail address:* `kapovich@math.uiuc.edu`

DEPARTMENT OF MATHEMATICS, UNIVERSITY OF ILLINOIS AT URBANA-CHAMPAIGN, 1409 WEST GREEN STREET, URBANA, IL 61801, USA

*E-mail address:* `schupp@math.uiuc.edu`